%% file: v-Dims.tex
\documentclass[12pt,reqno]{amsart}
\usepackage{fullpage,graphicx,dbnsymb,amsmath,amssymb,color}
\usepackage{hyperref}

\catcode`\@=11
\long\def\@makecaption#1#2{%
    \vskip 10pt
    \setbox\@tempboxa\hbox{
      \small\sf{\bfcaptionfont #1. }\ignorespaces #2}%
    \ifdim \wd\@tempboxa >\captionwidth {%
        \rightskip=\@captionmargin\leftskip=\@captionmargin
        \unhbox\@tempboxa\par}%
      \else
        \hbox to\hsize{\hfil\box\@tempboxa\hfil}%
    \fi}
\font\bfcaptionfont=cmssbx10 scaled \magstephalf
\newdimen\@captionmargin\@captionmargin=2\parindent
\newdimen\captionwidth\captionwidth=\hsize
\catcode`\@=12

\def\calAar{{\vec{\mathcal A}}}
\def\calP{{\mathcal P}}
\def\calV{{\mathcal V}}
\def\calW{{\mathcal W}}

\def\arXiv#1{{\href{http://front.math.ucdavis.edu/#1}{arXiv:#1}}}

\theoremstyle{plain}
\newtheorem{conjecture}{Conjecture}

\begin{document}
\newdimen\captionwidth\captionwidth=\hsize
\setcounter{secnumdepth}{4}

\title[Dimensions for Virtual Knots]{Some Dimensions of Spaces of Finite
Type Invariants of Virtual Knots}

\author{Dror~Bar-Natan}
\address{
  Department of Mathematics\\
  University of Toronto\\
  Toronto Ontario M5S 2E4\\
  Canada
}
\email{drorbn@math.toronto.edu}
\urladdr{http://www.math.toronto.edu/~drorbn}

\author{Iva Halacheva}
\address{
  Department of Mathematics\\
  University of Toronto\\
  Toronto Ontario M5S 2E4\\
  Canada
}
\email{iva.halacheva@utoronto.ca}

\author{Louis Leung}
\address{
  Department of Mathematics\\
  University of Toronto\\
  Toronto Ontario M5S 2E4\\
  Canada
}
\email{louis.leung@utoronto.ca}

\author{Fionntan Roukema}
\address{
  Department of Mathematics\\
  University of Pisa\\
  Largo Bruno Pontecorvo 5, 56127 Pisa\\
  Italy
}
\email{f.w.m.roukema.03@cantab.net}

\date{Sep.~28,~2009. 
  Electronic version and related files at~\cite{v-Dims}.
}


\begin{abstract}
  \input abstract.tex
\end{abstract}

\maketitle


\section{``Standard'' Virtual Knots} \label{sec:standard}

For ``classical'' finite type invariants of ordinary
knots, as defined by the schematic difference relation
{\large$\doublepoint\to\overcrossing-\undercrossing$} (see
e.g.~\cite{Bar-Natan:OnVassiliev}), it is well known that ``every weight
system integrates''. In other words, every linear functional on chord
diagrams which satisfies the 4T relation is the ``top derivative'' of some
finite type invariant. Indeed, this simple minded statement is the main
implication of the existence of the celebrated ``Kontsevich integral''
and of ``configuration space integrals'', and it is closely related to
``perturbative Chern-Simons theory'' and to the theory of ``Drinfel'd
associators'' (see overviews at~\cite{Bar-NatanStoimenow:Fundamental,
Bar-Natan:EMP}).

The purpose of this note is to support the conjecture that
the same is true in the context of ``v-knots'' or ``virtual
knots''~\cite{Kauffman:VirtualKnotTheory} (and in fact, also
in several closely related contexts). In this case finite type
invariants are defined by the schematic difference relation
{\large$\svslashoverback\to\slashoverback-\crossing$} (see
e.g.~\cite{GoussarovPolyakViro:VirtualKnots}).

We wrote a computer program (see~\cite{v-Dims}) to compute
the dimensions (``$\dim\calW_n$'') of spaces of weight systems
(of v-knots) of various degrees, and using the ``Polyak algebra''
of~\cite{GoussarovPolyakViro:VirtualKnots}, to also compute the dimensions
(``$\dim\calV_n/\calV_{n-1}$'', or more shortly, ``$\dim\calV_{n/n-1}$'')
of the spaces of finite type invariants (of v-knots) of various degrees
(modulo invariants of lower degree).  Here are the results:

\begin{center}
\begin{tabular}{||c||c|c|c|c|c|c||}
\multicolumn{7}{l}{Dimensions for round v-knots:} \\
\hline \hline
  $n$ & 0 & 1 & 2 & 3 & 4 & 5 \\
  \hline
  $\dim\calW_n$ & 1 & 0 & 0 & 1 & 4 & 17 \\
  \hline
  $\dim\calV_{n/n-1}$ & 1 & 0 & 0 & 1 & 4 & 17 \\
\hline\hline
\end{tabular}
\qquad
\begin{tabular}{||c||c|c|c|c|c|c||}
\multicolumn{7}{l}{Dimensions for long v-knots:} \\
\hline \hline
  $n$ & 0 & 1 & 2 & 3 & 4 & 5 \\
  \hline
  $\dim\calW_n$ & 1 & 0 & 2 & 7 & 42 & 246 \\
  \hline
  $\dim\calV_{n/n-1}$ & 1 & 0 & 2 & 7 & 42 & 246 \\
\hline\hline
\end{tabular}
\end{center}

\begin{conjecture} \label{conj:small}
The pattern of equalities appearing above continues. That is, every
weight system for v-knots comes from a finite type invariant of v-knots.
\end{conjecture}

\section{Variants}

\begin{figure}
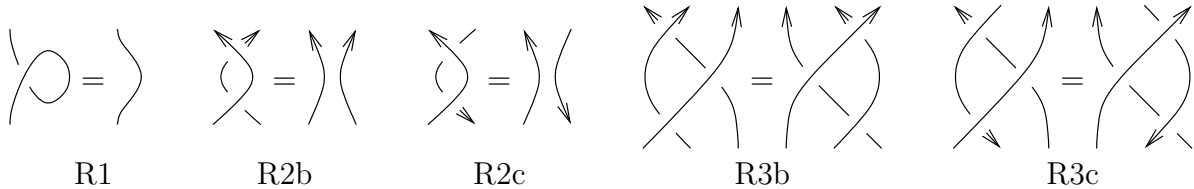

\[ \begin{array}{c}\input figs/RMoves.pstex_t \end{array} \]
\caption{Five types of Reidemeister moves. The ``b'' or ``braid-like''
moves R2b and R3b have all strands oriented the same way; such
configurations could
appear in a braid. The ``c'' or ``cyclic'' moves R2c and R3c contain a
planar domain whose boundary is oriented cyclically. Such configurations
cannot appear within a braid. Moves involving ``virtual crossings''
remain as in~\cite{Kauffman:VirtualKnotTheory} and are not shown here.}
\label{fig:RMoves}
\end{figure}

The theory of finite type invariants of ordinary knots is rather ``rigid''
--- it is the same for round or long knots, the framed and unframed
cases are not too different (in particular, a complete understanding of
one is equivalent to a complete understanding of the other), and there is
little else that can be tinkered with. This is not the case for virtual
knot theory --- round and long and framed and unframed appear to be quite
different, and there are several other ``parameters'' that can be turned on
and off at will, leading to a significant number of apparently different
``virtual knot theories''. In each such theory we start with a collection
of ``virtual knot diagrams'' and then mod it out by some ``Reidemeister
moves'' (see Figure~\ref{fig:RMoves}). Some of the possible choices follow:

\begin{itemize}
\item Skeleton choices: We can take the skeleton of our virtual knots to be
a circle (the ``round'' case) or a line (the ``long'' case). In the case of
a line, we may restrict our attention to virtual knot diagrams all of whose
(real) crossings are ``descending'' (a crossing is ``descending'' if the
first time it is visited along the parametrization of the knot it is
visited on the ``over'' strand).
\item R23 choices: We may mod out by all R2 and R3 moves (this is the
``standard'' case), or only by the
``braid-like'' moves R2b and R3b, or we may skip R3 moves altogether and
only mod out by R2b and R2c (``R2 only'').
\item R1 choices: We may or may not mod out by R1 moves.
\item Other choices: The ``Overcrossings Commute'' relation is studied
extensively in~\cite{Bar-Natan:WKO} and will not be studied here. ``Flat''
and ``free'' virtual knots are studied in~\cite{Manturov:FreeKnots,
Manturov:FreeKnotsLinks} and will not be studied here. ``Virtual braids''
are left for a future study.
\end{itemize}

Each such virtual knot theory has a notion of finite type invariants
(always defined by $\svslashoverback\to\slashoverback-\crossing$),
and each one has a notion of ``weight systems'' (see
Section~\ref{sec:ArrowDiagrams}). Hence the question ``does every weight
system come from a finite type invariant'' makes sense in many ways.
We have studied $18=3\times 3\times 2$ of these ways:

\begin{conjecture}[``18 in 1''] \label{conj:big}
For each skeleton choice (``round'', ``long'', or ``descending''), with
R2 and R3 given either the ``standard'' or the ``braid-like'' or the
``R2 only'' treatment, with or without R1, and for every natural number
$n$, every degree $n$ weight system comes from a type $n$ invariant.
\end{conjecture}

Using our program (see~\cite{v-Dims}) we have verified the above
conjecture for $n\leq 5$ in all 18 cases. Below we display the dimensions
of the spaces ``$\calV_{n/n-1}$'' of type $n$ invariants modulo invariants
of lower type (for each case). By our computer's hard work the dimensions
of the spaces ``$\calW_n$'' of weight systems are exactly the same,
so they do not require a separate table.  This equality of dimensions
for all $n$ is precisely the content of our conjecture. In all cases
$\dim\calV_0=\dim\calW_0=1$, so we only display the dimensions for
$n=1,2,3,4,5$:

\begin{center}\begin{tabular}{||l|l||c|c|c||}
\hline\hline
  \multicolumn{2}{||l||}{$\dim\calV_{n/n-1}$ / $\dim\calW_n$ for \ldots} &
  round v-knots &
  long v-knots &
  descending v-knots \\
\hline\hline
  standard R23 & mod R1 &
  $0,0,1,4,17^{\text{(\ref{com:standard})}}$ &
  $0,2,7,42,246^{\text{(\ref{com:standard})}}$ &
  $0,0,1,6,34$ \\
\cline{2-5}
  & no R1 &
  $1,1,2,7,29$ &
  $2,5,15,67,365$ &
  $1,1,2,8,42$ \\
\hline
  braid-like R23 & mod R1 &
  $0,0,1,4,17^{\text{(\ref{com:match})}}$ &
  $0,2,7,42,246^{\text{(\ref{com:match})}}$ &
  $0,0,1,6,34^{\text{(\ref{com:match})}}$ \\
\cline{2-5}
  & no R1 &
  $1,2,5,19,77$ &
  $2,7,27,139,813^{\text{(\ref{com:lie})}}$ &
  $1,2,6,24,120^{\text{(\ref{com:factorial})}}$ \\
\hline
  R2 only & mod R1 &
  $0,0,4,44,648$ &
  $0,2,28,420,7808$ &
  $0,0,2,18,174$ \\
\cline{2-5}
  & no R1 &
  $1,3,16,160,2248$ &
  $2,10,96,1332,23880$ &
  $1,2,9,63,570$ \\
\hline\hline
\end{tabular}\end{center}

\vskip 2mm
\par\noindent{\bf Comments. }
\begin{enumerate}
\item \label{com:standard} These are the ``standard'' virtual knots, as in
Section~\ref{sec:standard}.
\item \label{com:match} The equality of these numbers with the
numbers two rows above is a bit tricky. It is not true that R1
and the braid-like R23 imply the cyclic R23. Yet at the level of
arrow diagrams, FI and 6T do imply the XII relations (naming as in
Section~\ref{sec:ArrowDiagrams}). Thus the equality of $\dim\calW_n$'s is
obvious, and assuming Conjecture~\ref{conj:big} it implies the equality
of the $\dim\calV_{n/n-1}$'s.
\item \label{com:lie} The spaces measured in this box are dual to
(long arrow diagrams)/(6T relations), and these are the spaces most
closely related to Lie bi-algebras~\cite{Haviv:DiagrammaticAnalogue,
Leung:CombinatorialFormulas, Bar-Natan:WKO}. Thus in the long run this box
may prove to be the most important of the variants of ``virtual knots''
studied here.  \item \label{com:factorial} We can show that in this case
$\dim\calW_n\leq n!$ but we are missing the other inequality necessary
to prove that $\dim\calW_n=n!$.
\end{enumerate}

In our computations we used the Polyak algebra techniques
of~\cite{GoussarovPolyakViro:VirtualKnots} for the ``$\calV$'' spaces
and straightforward linear algebra for the ``$\calW$'' spaces described
below.  The typical $n=5$ computation involves determining the rank of
a very sparse matrix with a few tens of thousands rows and columns and
takes about an hour of computer time. The main part of the program was
written in Mathematica~\cite{Wolfram:Mathematica} with the heavier rank
computations delegated to LinBox~\cite{ProjectLinbox}.

\vskip 2mm \par\noindent{\bf Why bother?} Why bother with such an
``18 in 1'' conjecture?  We believe virtual knots in general,
and the question studied here on finite type invariants of
virtual knots in particular, might form the correct topological
framework for the study of quantum groups and the quantization of
Lie bi-algebras~\cite{Haviv:DiagrammaticAnalogue, Bar-Natan:WKO,
EtingofKazhdan:BialgebrasI}. But we are not sure yet which class of
virtual knots it is that we should study. Is it the standard class, as in
Section~\ref{sec:standard}, or is it the one closest to Lie bialgebras,
as in Comment~(\ref{com:lie}) above? Or maybe it is something else,
closely related?

Thus we believe that at least some of the 18 cases in
Conjecture~\ref{conj:big} are deeply interesting. As for the rest (the
cases involving ``R2 only'' or ``descending v-knots'', for example),
these may play two kinds of roles in the future:

\begin{enumerate}
\item The apparently harder cases, involving all Reidemeister moves and
round or long skeleta,  appear quite hard. The ``easier'' cases may
serve as ``baby versions'' that will force us to develop some of the
techniques which we may later use while studying the harder cases.
\item We certainly hope that eventually all 18 cases (and maybe a few more)
of Conjecture~\ref{conj:big} will find a uniform solution. Thus the
presence of so many variants of Conjecture~\ref{conj:big} may serve as a
further test of our understanding. Suppose we solved one of the ``harder''
cases. Is our solution modular enough to resolve all other cases as well?
\end{enumerate}

\section{Arrow Diagrams and Weight Systems} \label{sec:ArrowDiagrams}

This is a short descriptive section intended only to spell out in
brief, for reasons of completeness, the definitions of the spaces
$\calW_n$ of weight systems for each of the cases that we have
considered. The details of how and why the spaces described below
are related to finite type invariants of virtual knots can be found
in~\cite{GoussarovPolyakViro:VirtualKnots, Polyak:ArrowDiagrams}.

\begin{figure}
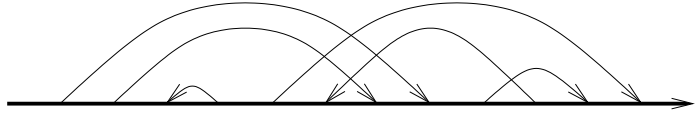

\[ \begin{array}{c}\input figs/ArrowDiagram.pstex_t \end{array} \]
\caption{
  A typical arrow diagram of degree 6 (meaning, having exactly 6
  arrows beyond the bolder ``skeleton'' line at the bottom).
} \label{fig:ArrowDiagram}
\end{figure}

\begin{figure}
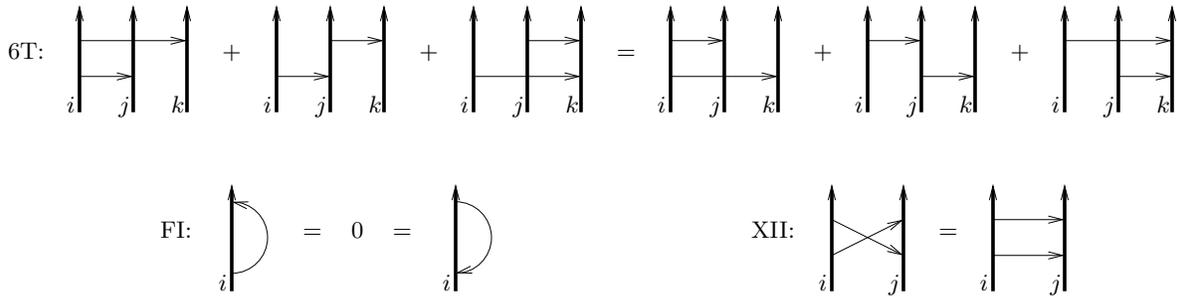

\[ \begin{array}{c}\input figs/ArrowRelations.pstex_t \end{array} \]
\caption{
  The 6T, FI and XII relations in standard ``skein'' notation --- only
  the varying parts of the diagrams involved are shown, their skeleton
  pieces (labeled $i$, $j$, and $k$ above) can be assembled along a
  long or a round skeleton in any way, and outside the parts shown,
  more arrows can be inserted.
} \label{fig:ArrowRelations}
\end{figure}

The spaces $\calW_n$ are always the duals $\calAar_n^\star$ of
spaces $\calAar_n$ of ``arrow diagrams'' modulo various kinds of
``arrow diagram relations''. A sample arrow diagram is shown in
Figure~\ref{fig:ArrowDiagram} and the templates for all the arrow diagram
relations that we consider are in Figure~\ref{fig:ArrowRelations}.

For long v-knots, $\calAar_n$ will consist of exactly the kind of arrow
diagrams shown in Figure~\ref{fig:ArrowDiagram}. For round v-knots we
replace the skeleton line by an oriented circle. For descending v-knots
the skeleton is again a long line, but for the diagrams in $\calAar_n$
we allow only those whose arrows are oriented the same way as their
skeleton (thus the sample diagram in Figure~\ref{fig:ArrowDiagram}
would be excluded because two of its arrows are oriented against the
orientation of the skeleton).

In the case of ``standard R23'', we impose both 6T and XII on
$\calAar_n$. In the case of ``braid-like R23'' we impose 6T but not
XII. In the ``R2 only'' case, we impose XII but not 6T.

We impose the FI relation in $\calAar_n$ iff we mod out by R1 at the level
of v-knots.

In the case of descending v-knots, we only impose 6T if $i<j<k$ (as sites
along the oriented skeleton), we only impose XII if $i<j$, and we only
impose the properly-oriented ``left half'' of FI.

The 6T and FI relations appear and are explained
in~\cite{GoussarovPolyakViro:VirtualKnots, Polyak:ArrowDiagrams}. For
all we know, this is relation XII's maiden appearance in the literature,
and thus an explanation is in order. Below are two brief derivations of
XII; the first direct and elementary, and the second using the Polyak
algebra. All relevant definitions are
in~\cite{GoussarovPolyakViro:VirtualKnots} and will not be repeated here.

\subsection{A Direct Derivation of XII}
The equality
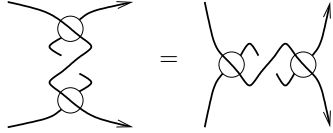
\[ \begin{array}{c}\input figs/Direct.pstex_t \end{array} \]
of semi-virtual tangles is easy to verify directly,
using the definitions of the semi-virtual crossing,
{\large$\svslashoverback=\slashoverback-\crossing$}, and using only
(virtual moves and) R2 moves (though both braid-like and cyclic ones).
But in arrow notation, this is exactly the XII relation.

\subsection{A Polyak Algebra Derivation of XII}

The Polyak algebra $\calP_n$ is defined
in~\cite{GoussarovPolyakViro:VirtualKnots}; it is a space of ``signed
arrow diagrams'' modulo relations that correspond to the Reidemeister
moves of knot theory. The relation corresponding to the R2 move is
\begin{equation} \label{eq:BarePolyakR2} \begin{array}{c}\input figs/BarePolyakR2.pstex_t \end{array} \end{equation}
Symbolically, with $a$ denoting the $+$ arrow and $b$ denoting the $-$
arrow, this is $ab+a+b=0$, or $b=-a-ab$. Solving for $b$ in terms of $a$
and remembering that in $\calP_n$ we mod out by degrees higher than $n$, we
get $b=-a+a^2-a^3+\cdots$ (a finite sum). Thus the negative arrows can be
eliminated in $\calP_n$ (this of course is very useful computationally, as
it lowers the number of arrow diagrams that one needs to consider by a
factor of about $2^n$).

But in Equation~\eqref{eq:BarePolyakR2} the orientation of the strands is
not specified, and indeed, for braid-like R2 moves these strands come
out with parallel orientations while for cyclic R2 moves they come out with
opposite orientations. Thus we get two different formulas for negative
arrows in terms of positive ones. The first, using parallel orientations
in~\eqref{eq:BarePolyakR2}, and dropping the signs from the positive
arrows, is
\begin{equation} \label{eq:FirstFormula} \begin{array}{c}\input figs/FirstFormula.pstex_t \end{array}
\end{equation}
In the second such formula, using opposite orientations 
in~\eqref{eq:BarePolyakR2}, we flip to the right the strand that was
oriented to the left at the cost of getting all the $a^k$ terms totally
twisted:
\begin{equation} \label{eq:SecondFormula} \begin{array}{c}\input figs/SecondFormula.pstex_t \end{array}
\end{equation}

Equating these two formulas and keeping only the lowest order terms that
don't cancel out, we get the XII relation:
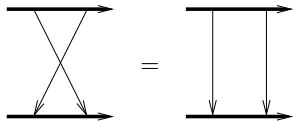
\[ \begin{array}{c}\input figs/XII.pstex_t \end{array} \]

The only benefit of the Polyak algebra derivation of XII is the following
comment: in the computation of $\dim\calP_n$ (which is the same as
$\dim\calV_n$) in the cases where all R2 moves are imposed, one may
restrict attention only to $+$ arrows, but then the full right-hand-sides
of~\eqref{eq:FirstFormula} and~\eqref{eq:SecondFormula} have to be set
equal, dropping only the terms of degree higher than $n$.

\section*{Disclaimer} Our computational results suggest what we believe
are interesting conjectures. Yet in programming, bugs are a fact of life.
An independent verification of our numbers, even without pushing beyond
degree 5, would lend further support to Conjectures~\ref{conj:small}
and~\ref{conj:big} and would be highly desirable.

\section*{Acknowledgement} We thank Bradford Hovinen, B.~David Saunders
and William Turner of Project LinBox for helping us with sparse matrix
computations, Karene Chu for some comments, and Carlo Petronio at the
University of Pisa for allowing Fionntan Roukema the freedom to partake
in the project. This work was partially supported by NSERC grant RGPIN
262178.

\end{document}

%% file: abstract.tex
We compute many dimensions of spaces of finite type invariants of virtual
knots (of several kinds) and the dimensions of the corresponding spaces of
``weight systems'', finding everything to be in agreement with the
conjecture that ``every weight system integrates''.

%% file: figs/RMoves.pstex_t
\begin{picture}(0,0)%
\includegraphics{figs/RMoves.pstex}%
\end{picture}%
%
%
\setlength{\unitlength}{3947sp}%
\begingroup\makeatletter\ifx\SetFigFont\undefined%
\gdef\SetFigFont#1#2#3#4#5{%
  \reset@font\fontsize{#1}{#2pt}%
  \fontfamily{#3}\fontseries{#4}\fontshape{#5}%
  \selectfont}%
\fi\endgroup%
\begin{picture}(7449,1137)(214,-286)
\put(6901,314){\makebox(0,0)[b]{\smash{{\SetFigFont{12}{14.4}{\rmdefault}{\mddefault}{\updefault}{\color[rgb]{0,0,0}$=$}%
}}}}
\put(751,-286){\makebox(0,0)[b]{\smash{{\SetFigFont{12}{14.4}{\rmdefault}{\mddefault}{\updefault}{\color[rgb]{0,0,0}R1}%
}}}}
\put(1951,-286){\makebox(0,0)[b]{\smash{{\SetFigFont{12}{14.4}{\rmdefault}{\mddefault}{\updefault}{\color[rgb]{0,0,0}R2b}%
}}}}
\put(3301,-286){\makebox(0,0)[b]{\smash{{\SetFigFont{12}{14.4}{\rmdefault}{\mddefault}{\updefault}{\color[rgb]{0,0,0}R2c}%
}}}}
\put(4951,-286){\makebox(0,0)[b]{\smash{{\SetFigFont{12}{14.4}{\rmdefault}{\mddefault}{\updefault}{\color[rgb]{0,0,0}R3b}%
}}}}
\put(6901,-286){\makebox(0,0)[b]{\smash{{\SetFigFont{12}{14.4}{\rmdefault}{\mddefault}{\updefault}{\color[rgb]{0,0,0}R3c}%
}}}}
\put(751,314){\makebox(0,0)[b]{\smash{{\SetFigFont{12}{14.4}{\rmdefault}{\mddefault}{\updefault}{\color[rgb]{0,0,0}$=$}%
}}}}
\put(1951,314){\makebox(0,0)[b]{\smash{{\SetFigFont{12}{14.4}{\rmdefault}{\mddefault}{\updefault}{\color[rgb]{0,0,0}$=$}%
}}}}
\put(3301,314){\makebox(0,0)[b]{\smash{{\SetFigFont{12}{14.4}{\rmdefault}{\mddefault}{\updefault}{\color[rgb]{0,0,0}$=$}%
}}}}
\put(4951,314){\makebox(0,0)[b]{\smash{{\SetFigFont{12}{14.4}{\rmdefault}{\mddefault}{\updefault}{\color[rgb]{0,0,0}$=$}%
}}}}
\end{picture}%

%% file: figs/ArrowDiagram.pstex_t
\begin{picture}(0,0)%
\includegraphics{figs/ArrowDiagram.pstex}%
\end{picture}%
%
%
\setlength{\unitlength}{3158sp}%
\begingroup\makeatletter\ifx\SetFigFont\undefined%
\gdef\SetFigFont#1#2#3#4#5{%
  \reset@font\fontsize{#1}{#2pt}%
  \fontfamily{#3}\fontseries{#4}\fontshape{#5}%
  \selectfont}%
\fi\endgroup%
\begin{picture}(5466,878)(-32,-106)
\end{picture}%

%% file: figs/ArrowRelations.pstex_t
\begin{picture}(0,0)%
\includegraphics{figs/ArrowRelations.pstex}%
\end{picture}%
%
%
\setlength{\unitlength}{2960sp}%
\begingroup\makeatletter\ifx\SetFigFontNFSS\undefined%
\gdef\SetFigFontNFSS#1#2#3#4#5{%
  \reset@font\fontsize{#1}{#2pt}%
  \fontfamily{#3}\fontseries{#4}\fontshape{#5}%
  \selectfont}%
\fi\endgroup%
\begin{picture}(9828,2511)(-464,-2839)
\put(2626,-1261){\makebox(0,0)[b]{\smash{{\SetFigFontNFSS{9}{10.8}{\rmdefault}{\mddefault}{\updefault}{\color[rgb]{0,0,0}$k$}%
}}}}
\put(2176,-1261){\makebox(0,0)[b]{\smash{{\SetFigFontNFSS{9}{10.8}{\rmdefault}{\mddefault}{\updefault}{\color[rgb]{0,0,0}$j$}%
}}}}
\put(1726,-1261){\makebox(0,0)[b]{\smash{{\SetFigFontNFSS{9}{10.8}{\rmdefault}{\mddefault}{\updefault}{\color[rgb]{0,0,0}$i$}%
}}}}
\put(4276,-1261){\makebox(0,0)[b]{\smash{{\SetFigFontNFSS{9}{10.8}{\rmdefault}{\mddefault}{\updefault}{\color[rgb]{0,0,0}$k$}%
}}}}
\put(3826,-1261){\makebox(0,0)[b]{\smash{{\SetFigFontNFSS{9}{10.8}{\rmdefault}{\mddefault}{\updefault}{\color[rgb]{0,0,0}$j$}%
}}}}
\put(3376,-1261){\makebox(0,0)[b]{\smash{{\SetFigFontNFSS{9}{10.8}{\rmdefault}{\mddefault}{\updefault}{\color[rgb]{0,0,0}$i$}%
}}}}
\put(5926,-1261){\makebox(0,0)[b]{\smash{{\SetFigFontNFSS{9}{10.8}{\rmdefault}{\mddefault}{\updefault}{\color[rgb]{0,0,0}$k$}%
}}}}
\put(5476,-1261){\makebox(0,0)[b]{\smash{{\SetFigFontNFSS{9}{10.8}{\rmdefault}{\mddefault}{\updefault}{\color[rgb]{0,0,0}$j$}%
}}}}
\put(5026,-1261){\makebox(0,0)[b]{\smash{{\SetFigFontNFSS{9}{10.8}{\rmdefault}{\mddefault}{\updefault}{\color[rgb]{0,0,0}$i$}%
}}}}
\put(7576,-1261){\makebox(0,0)[b]{\smash{{\SetFigFontNFSS{9}{10.8}{\rmdefault}{\mddefault}{\updefault}{\color[rgb]{0,0,0}$k$}%
}}}}
\put(7126,-1261){\makebox(0,0)[b]{\smash{{\SetFigFontNFSS{9}{10.8}{\rmdefault}{\mddefault}{\updefault}{\color[rgb]{0,0,0}$j$}%
}}}}
\put(6676,-1261){\makebox(0,0)[b]{\smash{{\SetFigFontNFSS{9}{10.8}{\rmdefault}{\mddefault}{\updefault}{\color[rgb]{0,0,0}$i$}%
}}}}
\put(9226,-1261){\makebox(0,0)[b]{\smash{{\SetFigFontNFSS{9}{10.8}{\rmdefault}{\mddefault}{\updefault}{\color[rgb]{0,0,0}$k$}%
}}}}
\put(8776,-1261){\makebox(0,0)[b]{\smash{{\SetFigFontNFSS{9}{10.8}{\rmdefault}{\mddefault}{\updefault}{\color[rgb]{0,0,0}$j$}%
}}}}
\put(8326,-1261){\makebox(0,0)[b]{\smash{{\SetFigFontNFSS{9}{10.8}{\rmdefault}{\mddefault}{\updefault}{\color[rgb]{0,0,0}$i$}%
}}}}
\put(5776,-2311){\makebox(0,0)[lb]{\smash{{\SetFigFontNFSS{9}{10.8}{\rmdefault}{\mddefault}{\updefault}{\color[rgb]{0,0,0}XII:}%
}}}}
\put(7426,-2311){\makebox(0,0)[b]{\smash{{\SetFigFontNFSS{9}{10.8}{\rmdefault}{\mddefault}{\updefault}{\color[rgb]{0,0,0}$=$}%
}}}}
\put(3226,-2761){\makebox(0,0)[b]{\smash{{\SetFigFontNFSS{9}{10.8}{\rmdefault}{\mddefault}{\updefault}{\color[rgb]{0,0,0}$i$}%
}}}}
\put(4726,-811){\makebox(0,0)[b]{\smash{{\SetFigFontNFSS{9}{10.8}{\rmdefault}{\mddefault}{\updefault}{\color[rgb]{0,0,0}$=$}%
}}}}
\put(1426,-811){\makebox(0,0)[b]{\smash{{\SetFigFontNFSS{9}{10.8}{\rmdefault}{\mddefault}{\updefault}{\color[rgb]{0,0,0}$+$}%
}}}}
\put(3076,-811){\makebox(0,0)[b]{\smash{{\SetFigFontNFSS{9}{10.8}{\rmdefault}{\mddefault}{\updefault}{\color[rgb]{0,0,0}$+$}%
}}}}
\put(6376,-811){\makebox(0,0)[b]{\smash{{\SetFigFontNFSS{9}{10.8}{\rmdefault}{\mddefault}{\updefault}{\color[rgb]{0,0,0}$+$}%
}}}}
\put(8026,-811){\makebox(0,0)[b]{\smash{{\SetFigFontNFSS{9}{10.8}{\rmdefault}{\mddefault}{\updefault}{\color[rgb]{0,0,0}$+$}%
}}}}
\put(-449,-811){\makebox(0,0)[lb]{\smash{{\SetFigFontNFSS{9}{10.8}{\rmdefault}{\mddefault}{\updefault}{\color[rgb]{0,0,0}6T:}%
}}}}
\put(976,-1261){\makebox(0,0)[b]{\smash{{\SetFigFontNFSS{9}{10.8}{\rmdefault}{\mddefault}{\updefault}{\color[rgb]{0,0,0}$k$}%
}}}}
\put(526,-1261){\makebox(0,0)[b]{\smash{{\SetFigFontNFSS{9}{10.8}{\rmdefault}{\mddefault}{\updefault}{\color[rgb]{0,0,0}$j$}%
}}}}
\put( 76,-1261){\makebox(0,0)[b]{\smash{{\SetFigFontNFSS{9}{10.8}{\rmdefault}{\mddefault}{\updefault}{\color[rgb]{0,0,0}$i$}%
}}}}
\put(1351,-2761){\makebox(0,0)[b]{\smash{{\SetFigFontNFSS{9}{10.8}{\rmdefault}{\mddefault}{\updefault}{\color[rgb]{0,0,0}$i$}%
}}}}
\put(6376,-2761){\makebox(0,0)[b]{\smash{{\SetFigFontNFSS{9}{10.8}{\rmdefault}{\mddefault}{\updefault}{\color[rgb]{0,0,0}$i$}%
}}}}
\put(7726,-2761){\makebox(0,0)[b]{\smash{{\SetFigFontNFSS{9}{10.8}{\rmdefault}{\mddefault}{\updefault}{\color[rgb]{0,0,0}$i$}%
}}}}
\put(8326,-2761){\makebox(0,0)[b]{\smash{{\SetFigFontNFSS{9}{10.8}{\rmdefault}{\mddefault}{\updefault}{\color[rgb]{0,0,0}$j$}%
}}}}
\put(6976,-2761){\makebox(0,0)[b]{\smash{{\SetFigFontNFSS{9}{10.8}{\rmdefault}{\mddefault}{\updefault}{\color[rgb]{0,0,0}$j$}%
}}}}
\put(826,-2311){\makebox(0,0)[lb]{\smash{{\SetFigFontNFSS{9}{10.8}{\rmdefault}{\mddefault}{\updefault}{\color[rgb]{0,0,0}FI:}%
}}}}
\put(2101,-2311){\makebox(0,0)[b]{\smash{{\SetFigFontNFSS{9}{10.8}{\rmdefault}{\mddefault}{\updefault}{\color[rgb]{0,0,0}$=$}%
}}}}
\put(2476,-2311){\makebox(0,0)[b]{\smash{{\SetFigFontNFSS{9}{10.8}{\rmdefault}{\mddefault}{\updefault}{\color[rgb]{0,0,0}$0$}%
}}}}
\put(2851,-2311){\makebox(0,0)[b]{\smash{{\SetFigFontNFSS{9}{10.8}{\rmdefault}{\mddefault}{\updefault}{\color[rgb]{0,0,0}$=$}%
}}}}
\end{picture}%

%% file: figs/Direct.pstex_t
\begin{picture}(0,0)%
\includegraphics{figs/Direct.pstex}%
\end{picture}%
%
%
\setlength{\unitlength}{2960sp}%
\begingroup\makeatletter\ifx\SetFigFont\undefined%
\gdef\SetFigFont#1#2#3#4#5{%
  \reset@font\fontsize{#1}{#2pt}%
  \fontfamily{#3}\fontseries{#4}\fontshape{#5}%
  \selectfont}%
\fi\endgroup%
\begin{picture}(2746,1097)(1704,-308)
\put(3076,239){\makebox(0,0)[b]{\smash{{\SetFigFont{9}{10.8}{\rmdefault}{\mddefault}{\updefault}{\color[rgb]{0,0,0}$=$}%
}}}}
\end{picture}%

%% file: figs/BarePolyakR2.pstex_t
\begin{picture}(0,0)%
\includegraphics{figs/BarePolyakR2.pstex}%
\end{picture}%
%
%
\setlength{\unitlength}{2960sp}%
\begingroup\makeatletter\ifx\SetFigFont\undefined%
\gdef\SetFigFont#1#2#3#4#5{%
  \reset@font\fontsize{#1}{#2pt}%
  \fontfamily{#3}\fontseries{#4}\fontshape{#5}%
  \selectfont}%
\fi\endgroup%
\begin{picture}(4677,966)(-32,-244)
\put(1201,164){\makebox(0,0)[b]{\smash{{\SetFigFont{9}{10.8}{\rmdefault}{\mddefault}{\updefault}{\color[rgb]{0,0,0}$+$}%
}}}}
\put(2701,164){\makebox(0,0)[b]{\smash{{\SetFigFont{9}{10.8}{\rmdefault}{\mddefault}{\updefault}{\color[rgb]{0,0,0}$+$}%
}}}}
\put(4276,164){\makebox(0,0)[lb]{\smash{{\SetFigFont{9}{10.8}{\rmdefault}{\mddefault}{\updefault}{\color[rgb]{0,0,0}$=0$}%
}}}}
\put(301,464){\makebox(0,0)[b]{\smash{{\SetFigFont{9}{10.8}{\rmdefault}{\mddefault}{\updefault}{\color[rgb]{0,0,0}$+$}%
}}}}
\put(751,464){\makebox(0,0)[b]{\smash{{\SetFigFont{9}{10.8}{\rmdefault}{\mddefault}{\updefault}{\color[rgb]{0,0,0}$-$}%
}}}}
\put(2026,464){\makebox(0,0)[b]{\smash{{\SetFigFont{9}{10.8}{\rmdefault}{\mddefault}{\updefault}{\color[rgb]{0,0,0}$+$}%
}}}}
\put(3526,464){\makebox(0,0)[b]{\smash{{\SetFigFont{9}{10.8}{\rmdefault}{\mddefault}{\updefault}{\color[rgb]{0,0,0}$-$}%
}}}}
\end{picture}%

%% file: figs/FirstFormula.pstex_t
\begin{picture}(0,0)%
\includegraphics{figs/FirstFormula.pstex}%
\end{picture}%
%
%
\setlength{\unitlength}{2960sp}%
\begingroup\makeatletter\ifx\SetFigFont\undefined%
\gdef\SetFigFont#1#2#3#4#5{%
  \reset@font\fontsize{#1}{#2pt}%
  \fontfamily{#3}\fontseries{#4}\fontshape{#5}%
  \selectfont}%
\fi\endgroup%
\begin{picture}(5583,1026)(-182,-1624)
\put(4051,-1186){\makebox(0,0)[b]{\smash{{\SetFigFont{9}{10.8}{\rmdefault}{\mddefault}{\updefault}{\color[rgb]{0,0,0}$-$}%
}}}}
\put(1051,-1186){\makebox(0,0)[b]{\smash{{\SetFigFont{9}{10.8}{\rmdefault}{\mddefault}{\updefault}{\color[rgb]{0,0,0}$=$}%
}}}}
\put(2701,-1186){\makebox(0,0)[b]{\smash{{\SetFigFont{9}{10.8}{\rmdefault}{\mddefault}{\updefault}{\color[rgb]{0,0,0}$+$}%
}}}}
\put(1351,-1186){\makebox(0,0)[b]{\smash{{\SetFigFont{9}{10.8}{\rmdefault}{\mddefault}{\updefault}{\color[rgb]{0,0,0}$-$}%
}}}}
\put(5401,-1186){\makebox(0,0)[lb]{\smash{{\SetFigFont{9}{10.8}{\rmdefault}{\mddefault}{\updefault}{\color[rgb]{0,0,0}$+\cdots$}%
}}}}
\put(376,-886){\makebox(0,0)[b]{\smash{{\SetFigFont{9}{10.8}{\rmdefault}{\mddefault}{\updefault}{\color[rgb]{0,0,0}$-$}%
}}}}
\end{picture}%

%% file: figs/SecondFormula.pstex_t
\begin{picture}(0,0)%
\includegraphics{figs/SecondFormula.pstex}%
\end{picture}%
%
%
\setlength{\unitlength}{2960sp}%
\begingroup\makeatletter\ifx\SetFigFont\undefined%
\gdef\SetFigFont#1#2#3#4#5{%
  \reset@font\fontsize{#1}{#2pt}%
  \fontfamily{#3}\fontseries{#4}\fontshape{#5}%
  \selectfont}%
\fi\endgroup%
\begin{picture}(5583,1026)(-182,-1624)
\put(4051,-1186){\makebox(0,0)[b]{\smash{{\SetFigFont{9}{10.8}{\rmdefault}{\mddefault}{\updefault}{\color[rgb]{0,0,0}$-$}%
}}}}
\put(1051,-1186){\makebox(0,0)[b]{\smash{{\SetFigFont{9}{10.8}{\rmdefault}{\mddefault}{\updefault}{\color[rgb]{0,0,0}$=$}%
}}}}
\put(2701,-1186){\makebox(0,0)[b]{\smash{{\SetFigFont{9}{10.8}{\rmdefault}{\mddefault}{\updefault}{\color[rgb]{0,0,0}$+$}%
}}}}
\put(1351,-1186){\makebox(0,0)[b]{\smash{{\SetFigFont{9}{10.8}{\rmdefault}{\mddefault}{\updefault}{\color[rgb]{0,0,0}$-$}%
}}}}
\put(5401,-1186){\makebox(0,0)[lb]{\smash{{\SetFigFont{9}{10.8}{\rmdefault}{\mddefault}{\updefault}{\color[rgb]{0,0,0}$+\cdots$}%
}}}}
\put(376,-886){\makebox(0,0)[b]{\smash{{\SetFigFont{9}{10.8}{\rmdefault}{\mddefault}{\updefault}{\color[rgb]{0,0,0}$-$}%
}}}}
\end{picture}%

%% file: figs/XII.pstex_t
\begin{picture}(0,0)%
\includegraphics{figs/XII.pstex}%
\end{picture}%
%
%
\setlength{\unitlength}{2960sp}%
\begingroup\makeatletter\ifx\SetFigFont\undefined%
\gdef\SetFigFont#1#2#3#4#5{%
  \reset@font\fontsize{#1}{#2pt}%
  \fontfamily{#3}\fontseries{#4}\fontshape{#5}%
  \selectfont}%
\fi\endgroup%
\begin{picture}(2466,1026)(2818,-1624)
\put(4051,-1186){\makebox(0,0)[b]{\smash{{\SetFigFont{9}{10.8}{\rmdefault}{\mddefault}{\updefault}{\color[rgb]{0,0,0}$=$}%
}}}}
\end{picture}%